\documentclass{amsart}     
\usepackage{amsmath,amssymb,amsfonts} 
\usepackage[dvips]{graphicx}
\usepackage{pst-all}
\usepackage{graphics}                 
\usepackage{color}
\newpsobject{showgrid}{psgrid}{subgriddiv=1,griddots=5,gridlabels=0pt}
\usepackage{amssymb, amsmath, amscd}

\newtheorem{Thm}{Theorem} 
\newtheorem{prop}[Thm]{Proposition} 
\newtheorem{lem}[Thm]{Lemma} 
\newtheorem{cor}[Thm]{Corollary} 
\newtheorem{defi}[Thm]{Definition}

\newcommand{\Z}{\mathbb{Z}}

\newcommand{\R}{\mathbb{R}}

\begin{document}

\author{Paolo Salvatore and Roberto Tauraso}
\title{The operad Lie is free}
\date{20-02-2008}
\maketitle

\begin{abstract}
We show that the operad $Lie$ is free as a non-symmetric operad.
Then we study the generating series counting the operadic generators, 
finding a recursive formula for its coefficients, and 
showing that the asymptotic density of the operadic generators is 
$1/e$.
\end{abstract}

\section{Introduction}

The operad $Lie$ is the symmetric operad encoding the structure of a Lie algebra. It has an antisymmetric binary generator, the bracket, and a ternary  relation, the Jacobi identity.
For standard facts and notations on operads we refer to \cite{mss}.
Over a commutative ring $R$,
the $k$-term of the operad $Lie(k)$ is the subgroup of the free Lie algebra on $k$
generators $x_1,\dots,x_k$ spanned by the words where each generator appears exactly once. 
The action of the symmetric groups on $k$ letters exchanges the indices.
It is well known that $Lie(k)$ is a free $R$-modules of rank $(k-1)!$
\cite{Sinha}
The algebras over $Lie$ in the operadic sense are exactly Lie algebras in the usual sense, unless $R$ has 2-torsion.

A topological interpretation of $Lie$ was discovered by Fred Cohen.
Let $D_n$ be the little $n$-discs operad. The top homology
group of $D_n(k)$ is $$H_{(k-1)(n-1)}(D_n(k),R) \cong R^{(k-1)!} $$
Cohen proved that for $n>2$ odd the induced operad structure on
the top homology groups of $D_n$ is isomorphic to $Lie$.
This story is well explained by D. Sinha in \cite{Sinha}.

If we forget about the action of the symmetric groups, we can regard 
$Lie$ as a non-symmetric operad. 
A motivation to study $Lie$ in this sense comes from knot theory. Lambrechts Turchin and Volic proved in \cite{LTV}
that the rational homology of the space of long knots in $\R^n$ is the value on the operad $H_*(D_n)$ of a functor (Hochschild homology) defined naturally on non-symmetric operads with multiplication. The main result of this note is the following.  

\begin{Thm}
The operad {\rm Lie} is a free non-symmetric operad.
\end{Thm}

Section 2 is devoted to the proof of the theorem.
We start by constructing an operad $L$ in the category of sets 
that spans linearly $Lie$. Its elements are those iterated brackets
in $Lie$ such that the smallest index and the largest index inside
each bracket lie respectively on the left and on the right hand side. 
We recall an explicit construction of the free non-symmetric operad.
Then we show that the
operad $L$ is free, generated by 
a sequence of sets $P=(P(n))$. The elements of $P$, that we call
{\em primes}, are those iterated brackets such that the indices
inside each (non outer) bracket do {\em not} form an interval of consecutive integers.

In section 3 we study the sequence of the cardinalities
$b_n=|P(n)|$ of the sets of prime generators.
We show in Theorem \ref{formula} that the following recursive formula holds:

$$b_2=1,\; b_n=\sum_{k=2}^{n-2} ((k+1)b_{k+1}+b_{k})b_{n-k}  \quad\mbox{for $n\geq 3$}.$$

The formula suggests that there might
be an explicit inductive way 
to construct the prime generators.

We also study the asymptotic density of the prime generators
showing in Theorem \ref{density} that

$${|P(n)|\over |L(n)|}={b_n\over (n-1)!}={1\over e}\left(1-{3\over n}-{5\over 2n^2}+O\left({1\over n^3}\right)\right).$$

\

The counting sequence of our prime generators is closely related to the counting sequence of 
the SIF permutations studied
by Callan in \cite{Callan}. We derive in Corollary \ref{sif} a formula for the asymptotic density of the SIF permutations, that 
was conjectured in \cite{Callan}.

\section{Proof of the main theorem}

We will define, as stated in the introduction,
an operad $L$ in the category of sets.

\begin{defi}
Let $L(k)$ be the set of formal expressions obtained by applying iterated binary brackets to $k$ symbols $x_1,\dots,x_k$, such that

1) Each symbol appears exactly once 

2) The smallest index inside a bracket appears on the left hand side, and 
the largest index appears on the right hand side.
\end{defi}

For example the expression $[x_1,[x_2,x_3]]$ is in $L(3)$, but 
$[x_2,[x_1,x_3]]$ is not in $L(3)$ because $x_1$ is not on the left hand
side of the outer bracket.

By induction each expression in $k$ symbols must involve exactly 
$k-1$ brackets.
For example $L(1)$ contains only the unbracketed expression $x_1$,
$L(2)$ contains only the bracket $[x_1,x_2]$, and $L(3)$ contains the elements $[x_1,[x_2,x_3]], \, [[x_1,x_2],x_3]$.

\begin{prop}
The set $L(k)$ has cardinality $(k-1)!$
\end{prop}
\begin{proof}
By induction on $k$. This is true for $k=2$.
Suppose that this is true for $k<n$. An expression $A \in L(n)$
has the form $A=[A_1,A_2]$, where $A_1$ is an expression involving
symbols with indexes $1=a_1<\dots<a_j$ and $A_2$ involves symbols 
with indexes $b_1<\dots<b_{n-j}=n$, for some 
$1 \leq j \leq n-1$. If we replace each $a_i$ by $i$ in the expression $A_1$ we get an element $A'_1 \in L(j)$, and similarly replacing each $b_i$ by $i$ in $A_2$ gives $A'_2 \in L(n-j)$. 

\noindent How many expressions do we obtain in this way in $L(n)$ for fixed $j,A'_1$ and $A'_2$?
There are $\binom{n-2}{j-1}$ ways of shuffling the $a_i's$ and the $b_i's$
into two disjoint sets of cardinality $j$ and $n-j$,
with 1 belonging to the first and $n$
to the second. Then, by inductive hypothesis, if we fix $j$ and let $A'_1,A'_2$ vary we have $$\binom{n-2}{j-1}(j-1)!(n-j-1)!=(n-2)!$$
expressions, and summing over $j$ we get $(n-1)(n-2)!=(n-1)!$ expressions. 
\end{proof}

The collection $L(k)$ forms a non-symmetric operad in the category of sets, where the composition operation $A \circ_i B$ replaces the variable $x_i$ in $A$ by the formal expression $B$, with its indices shifted by $i-1$, and shifts by $b-1$ the indices of $A$ larger than $i$, where $b$ is the number of symbols in $B$.

For example $$[[x_1,x_3],[x_2,x_4]] \circ_3 [x_1,x_2] =
[[x_1,[x_3,x_4]],[x_2,x_5] ]$$

The unit of the operad is
the expression $x_1$.

\
We recall the definition of a free non-symmetric operad in the category of sets and $R$-modules. It will be sufficient for our purpose to consider reduced operads with no generators in degree (or arity) 0 and 1.

\begin{defi}
A reduced tree $T$ on $k$ leaves is a collection of subsets of $\{x_1,\dots,x_k \}$, the {\em vertices}, such that 

\begin{enumerate}
\item The indexes of the symbols of a vertex form a sequence of consecutive numbers

$\{i+1,\dots,i+k\}$ of cardinality $k \geq 2$;
\item The full set $\{x_1,\dots,x_k \}$ is a vertex;

\item Any symbol $x_i$ belongs to a vertex of $T$;

\item Given two distinct vertices $v,w$ of $T$, either $v \cap w = \emptyset$, or $v \subset w$, or $w \subset v$.
\end{enumerate}
\end{defi}

The {\em valence} $|v|$ of a vertex $v$ is its cardinality.
The terminology is justified because we can associate to each tree a planar directed graph that is a tree in the sense that it has no cycles, 
and each vertex in our sense yields a vertex of that graph,
with a number of incoming edges equal to the valence.
An example is given in the figure.

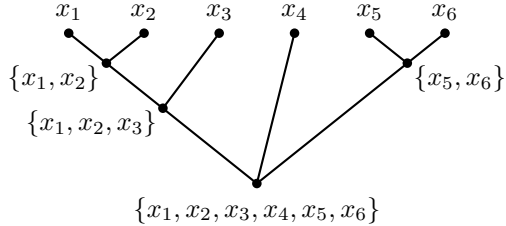
\begin{center}
\begin{pspicture}[.5](0,-0.7)(5,2.7)
\psline(2.5,0)(0,2)
\psline(2.5,0)(5,2)
\psline(2.5,0)(3,2)
\psline(0.5,1.6)(1,2)
\psline(4.5,1.6)(4,2)
\psline(1.25,1)(2,2)
\psdots(2.5,0)(0,2)(1,2)(2,2)(3,2)(4,2)(5,2)(0.5,1.6)(4.5,1.6)(1.25,1)
\uput[u](0,2){$x_1$}
\uput[u](1,2){$x_2$}
\uput[u](2,2){$x_3$}
\uput[u](3,2){$x_4$}
\uput[u](4,2){$x_5$}
\uput[u](5,2){$x_6$}
\rput(-0.2,1.4){$\{x_1,x_2\}$}
\rput(5.2,1.4){$\{x_5,x_6\}$}
\rput(0.3,0.8){$\{x_1,x_2,x_3\}$}
\uput[d](2.5,0){$\{x_1,x_2,x_3,x_4,x_5,x_6\}$}
\end{pspicture}
\end{center}

Let $M_k$ be the set of all reduced trees on $k$ leaves.
\begin{defi}
For given trees $T_1 \in M_k, \,T_2 \in M_l$ and $1 \leq i \leq k$, we define a new tree
$T_1 \circ_i T_2 \in M_{k+l-1}$ with the following vertices:

1) For each vertex $v \in T_1$ a vertex $v' \in T_1 \circ_i T_2$ containing 

-the symbol $x_j$ if $x_j \in v$ and $j<i$

-the symbol $x_{j+l-1}$ if $x_j \in v$ and $j>i$

-the symbols $x_{j+i-1}$ for $j=1,\dots,l$ if $x_i \in v$

2) For each vertex $u \in T_2$ a vertex $u''\in T_1 \circ_i T_2$ containing

-the symbol $x_{j+i-1}$ if $x_j \in u$.

\end{defi}

This gives a natural bijection of collections $T_1 \coprod T_2 \cong T_1 \circ_i T_2$.

Given a sequence of sets 
$X=(X_n)_{n>1}$, the free operad 
$F(X)$ is defined in degree (or arity) $k>1$ by
$$F(X)(k)=\coprod_{T \in M_k} \prod_{v \in T}X_{|v|}$$
We call the summand indexed by a tree $T \in M_k$ the {\em stratum}
$F(X)_T$.

The $\circ_i$ composition $F(X)(k) \times F(X)(l) \to F(X)(k+l-1)$ is
the inclusion that
identifies the product $F(X)_{T_1} \times F(X)_{T_2}$
of the strata 
 indexed respectively by $T_1 \in M_k$ and $T_2 \in M_l$ to the stratum 
$F(X)_{T_1 \circ_i T_2}$
indexed by $T_1 \circ_i T_2
\in M_{k+l-1}$.

The construction of a free operad in $R$-modules is entirely similar,
except that the disjoint union and the cartesian product of sets are replaced respectively by the direct sum 
and the tensor product of $R$-modules .

In particular the free $R$-module functor $X \mapsto R[X]$ sending sequences of sets to sequences of $R$-modules commutes with 
the free operad construction. 

\begin{prop} \label{pro}
The operad $L$ is free (up to isomorphism).
\end{prop}
\begin{proof}
Let us be given a formal expression $A \in L(k)$. We say that 
a bracket $b$ of $A$ is {\em connected} if the
set of indices inside $b$ is an interval of consecutive integers.
The indices do not need to appear in increasing order from left to right in the bracket. 
We say that a formal expression $A$ is {\em prime} if the outer bracket is the only connected bracket it contains.
For example $$[[x_1,x_3],[[x_2,x_4],x_5]] $$ is prime but 
\begin{equation} \label{deco}
[[[x_1,x_3],[x_2,x_4]],x_5] 
\end{equation}
 is not prime because the bracket
$[[x_1,x_3],[x_2,x_4]]$ is connected.

One can associate to each element of $L(k)$ a chord diagram, drawing
the indices 1,\dots,k on a line, and drawing 
for each bracket a chord from the smallest index inside the bracket to the largest index inside the bracket. 

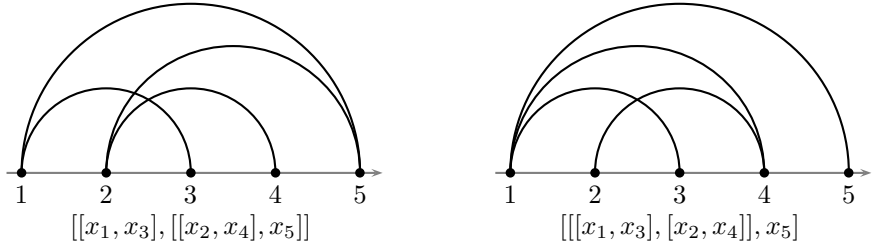
\begin{center}
\begin{pspicture}[.5](0,-1)(11.5,3)

\psarc(2.25,0){2.25}{0}{180}
\psarc(1.125,0){1.125}{0}{180}
\psarc(2.25,0){1.125}{0}{180}
\psarc(2.8125,0){1.6875}{0}{180}
\psline[linecolor=gray]{->}(-0.2,0)(4.8,0)
\psdots(0,0)(1.125,0)(2.25,0)(3.375,0)(4.5,0)
\uput[d](0,0){$1$}
\uput[d](1.125,0){$2$}
\uput[d](2.25,0){$3$}
\uput[d](3.375,0){$4$}
\uput[d](4.5,0){$5$}
\uput[d](2.25,-0.375){$[[x_1,x_3],[[x_2,x_4],x_5]]$}

\rput(6.5,0){
\psarc(2.25,0){2.25}{0}{180}
\psarc(1.125,0){1.125}{0}{180}
\psarc(2.25,0){1.125}{0}{180}
\psarc(1.6875,0){1.6875}{0}{180}
\psline[linecolor=gray]{->}(-0.2,0)(4.8,0)
\psdots(0,0)(1.125,0)(2.25,0)(3.375,0)(4.5,0)
\uput[d](0,0){$1$}
\uput[d](1.125,0){$2$}
\uput[d](2.25,0){$3$}
\uput[d](3.375,0){$4$}
\uput[d](4.5,0){$5$}
\uput[d](2.25,-0.375){$[[[x_1,x_3],[x_2,x_4]],x_5]$}
}
\end{pspicture}
\end{center}

 Let $P(n) \subset L(n)$ be the set of all prime expressions in $n$
 symbols, and let us consider the sequence of sets 
 $P=(P(n))$. We claim that $F(P)$ is
isomorphic to
$L$. The universal property defines an operad map $\theta: F(P) \to L$ out of the
inclusion $P \subset L$. We show that $\theta$ is an isomorphism constructing its inverse $\psi$.
Given $A \in L(k)$,
let $T_A \in M_k$ be the tree that has exactly a vertex $v_b=\{x_i,\dots,x_j\}$ for each connected bracket
$b$ of $A$ with set of indices $\{i,\dots,j\}$. 

To each such vertex we associate a prime expression $p_b$ obtained as follows: let $k'$ be the number of maximal connected brackets contained properly in $b$. Let $k''$ be the number of symbols in $b$ that are not contained in those $k'$ maximal brackets. Then there is a unique monotone surjective map $\phi:\{i,\dots,j\} \to \{1,\dots,k'+k''\}$ that is constant on the set of indices of each maximal bracket. 
The prime expression $p_b$ is obtained from $b$ replacing each index $t$ outside a maximal bracket by $\phi(t)$, and each maximal bracket $b'$ by a single symbol indexed by $\phi(u)$, where $u$ is any index inside $b'$. 
For example the expression \eqref{deco} corresponds to the tree
with a vertex $v=\{x_1,x_2,x_3,x_4,x_5\}$ labelled by $[x_1,x_2]$ and
a vertex $v'=\{x_1,x_2,x_3,x_4\}$ labelled by $[[x_1,x_3],[x_2,x_4]]$.

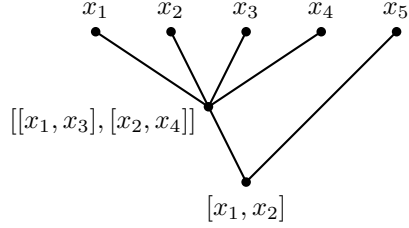
\begin{center}
\begin{pspicture}[.5](0,0.5)(4,3.7)
\psline(1.5,2)(0,3)
\psline(1.5,2)(1,3)
\psline(1.5,2)(2,3)
\psline(1.5,2)(3,3)
\psline(1.5,2)(2,1)
\psline(2,1)(4,3)
\psdots(0,3)(1,3)(2,3)(3,3)(4,3)(1.5,2)(2,1)
\uput[u](0,3){$x_1$}
\uput[u](1,3){$x_2$}
\uput[u](2,3){$x_3$}
\uput[u](3,3){$x_4$}
\uput[u](4,3){$x_5$}
\rput(0.1,1.8){$[[x_1,x_3],[x_2,x_4]]$}
\uput[d](2,1){$[x_1,x_2]$}
\end{pspicture}
\end{center}

The collection $\{p_b\}_{b \in T(A)}$ defines an element
 $\psi_k(A) \in F(P)(k)$. By construction $\psi=(\psi_k)$ is exactly the inverse of $\theta$.
\end{proof}

There is a tautological map of operads of $R$-modules 
$\beta: R[L] \to Lie$, sending a formal expression to the same expression inside $Lie$.

\begin{lem} \label{le}
The map $\beta$ is an isomorphism. 
\end{lem}

\begin{proof}
It is sufficient to establish the lemma for $R=\Z$.
We prove first the surjectivity.
We say that a bracket of an expression of iterated brackets has concentric index $i$ if it is properly contained in exactly $i-1$ brackets. For example the outer bracket has concentric index 1.  
Let us filter $Lie(k)$ by submodules $Lie_i(k)$ so that
$Lie_i(k)$ is linearly generated by iterated brackets where
in each bracket of index $\leq i$ the smallest index appears on the left and the largest on the right. In particular $Lie_0(k)=Lie(k), \,
Lie_{i+1}(k) \subseteq Lie_i(k)$ and $Lie_{k-1}(k)=Im(\beta_k)$.
We prove first that $Lie(k)=Lie_1(k)$.
Let $A \in Lie(k)$ be a linear generator, a combination of $k$ symbols using $k-1$ brackets.
We say that $A$ has weight $i$ if the smallest index 1 and the largest index $k$ share exactly $i$ (concentric) brackets.
For $i>1$ an application of the Jacobi identity shows that $A$ is the linear combination of two elements of weight $i-1$. 
For example $[[x_1,x_3],x_2]=-[x_3,[x_2,x_1]]-[x_2,[x_1,x_3]]$ reduces
an element of weight 1 as a combination of two elements of weight 0.
By induction $A$ is a linear combination of elements of weight 0, and 
by antisymmetry it is a linear combination of elements with the smallest index 1 on the left hand side of the outer bracket, and the largest index $k$ on the right hand side of the outer bracket. 
Thus $Lie(k)=Lie_1(k)$.

The same principle applied to brackets of concentric index i shows that 
$Lie_i(k)=Lie_{i+1}(k)$ for $i=1,\dots,k-2$.

Thus $Lie(k)=Lie_{k-1}(k)=Im(\beta_k)$ and $\beta_k$ is surjective. 
But the domain and the range of $\beta_k$ have both rank $(k-1)!$ in 
degree $k$ and so $\beta_k$ is an isomorphism for any $k$.
\end{proof} 

\

Theorem 1 follows from Proposition \ref{pro} and Lemma \ref{le}.
\newpage
\section{Counting the prime generators}

In this section we study the formal power series counting the number of prime generators.

\begin{lem} 
Let $(X_n)$ be a sequence of finite sets of cardinality $|X_n|=\beta_n$, for $n \geq 2$.
Let $FX$ be the free (reduced) operad generated by the sequence $X=(X_n)$. If $\alpha_n=|FX(n)|$ for $n \geq 1$, then the formal power
series $\alpha(x)=\sum_{n=1}^{\infty} \alpha_n x^n$ and 
$\beta(x)=\sum_{n=2}^{\infty} \beta_n x^n$ satisfy the identity
$$\beta(\alpha(x))+x=\alpha(x)$$
\end{lem}
 
\begin{proof}
Any element of $FX$ is either the unit $e \in FX(1)$, or it can be written uniquely as operadic composition $x(y_1,\dots,y_m)$, for some
$m \geq 2, \, x \in X_m$ and $y_i \in FX$, with $i=1,\dots,m$.
\end{proof} 
Actually the lemma is a special case of the non-symmetric version of
a result by Ginzburg and Kapranov, Theorem 3.3.2 in \cite{GK}, relating the Poincar\'e series of a differential graded operad to the series of its dual, or bar construction.

\

We know that the formal power series of $Lie$ is 
$$F(x)=\sum_{n=1}^{\infty}(n-1)! x^n $$
Let $b_n$ be the number of prime expressions in $L(n)$ and let
$$B(x)=-x+\sum_{n=2}^{\infty}b_n x^n.$$
The lemma applied to $\alpha(x)=F(x)$ and $\beta(x)=B(x)+x$ yields
the identity $-B(F(x))=x$ that is $B(x)=-F^{<-1>}(x)$, where the formal inverse 
$F^{<-1>}$ of $F$ exists because $F(x)=x+$  (higher order terms).

\begin{Thm} \label{formula}
The formal power series $B(x)$ counting the operadic generators 
of {\rm Lie} satisfies the differential equation
$$xB'(x)+(B'(x)+B(x))B(x)=0$$
and the following recursive formula holds:
$$b_2=1,\; b_n=\sum_{k=2}^{n-2} ((k+1)b_{k+1}+b_{k})b_{n-k}  \quad\mbox{for $n\geq 3$}.$$
\end{Thm}
\begin{proof}
By differentiating the identity $-B(F(x))=x$ we obtain
$$-B'(F(x))F'(x)=1.$$
Since $F'(x)=(F(x)-x)/x^2$ then
$$-B'(F(x))(F(x)-x)=x^2$$
and, by using $F(-B(x))=x$, we get 
$$xB'(x)+(B'(x)+B(x))B(x)=0.$$
Interpreting this identity termwise yields
$$[x^n]\left(xB'(x)+(B'(x)+B(x))B(x)\right)=nb_n+\sum_{k=0}^n ((k+1)b_{k+1}+b_k)b_{n-k}=0$$
and, since $b_0=0$, $b_1=-1$ and $b_2=1$, we have that
$$\sum_{k=2}^{n-2} ((k+1)b_{k+1}+b_k)b_{n-k}=-nb_n-(-1)b_{n}-(2-1)b_{n-1}-(nb_{n}+b_{n-1})(-1)=b_n.$$
\end{proof}

We obtain thus 
$$\sum_{n=2}^{\infty}b_n x^n=x^2+x^4+4x^5+22x^6+144x^7+1089x^8+9308x^9+88562x^{10}+\cdots$$
For more terms of $\{b_n\}_{n\geq 2}$ see the sequence A134988 in Sloane's OEIS \cite{Sloane}.

\begin{Thm} \label{density}
The asymptotics density of the operadic generators of ${\rm Lie}$ is
$${|P(n)|\over|dim({\rm Lie}(n))|}={b_n\over (n-1)!}=e^{-1}\left(1-{3\over n}-{5\over 2n^2}+O\left({1\over n^3}\right)\right).$$
\end{Thm}
\begin{proof}
Since $B(x)=-F^{<-1>}(x)$, then by the Lagrange inversion formula
\begin{eqnarray*}
{b_n\over (n-1)!}&=&-{1\over (n-1)!}[x^n]F^{<-1>}(x)=-{1\over n!}[x^{n-1}]\left({x\over F(x)}\right)^n\\
&=&-{1\over n!}[x^{n-1}]\left({1+A(x)}\right)^{-(n-1)-1}=-{1\over n!}(\alpha e^{\alpha a_1 \gamma}(n-1)a_{n-1}+O(a_{n-1}))\\
&=&-{1\over n!}(-e^{-1}(n-1)(n-1)!+O((n-1)!))=e^{-1}+O(1/n)
\end{eqnarray*}
where we applied Theorem 1(i) of \cite{Bender} with $A(x)=\sum_{k=1}^{\infty}k!x^k$, $\alpha=-1$, $\beta=-1$, $\gamma=1$.
Let $p_n=b_n/(n-1)!$ and 
$$p_n=e^{-1}\left(1+{c_1\over n}+{c_2\over n^2}+O\left({1\over n^3}\right)\right).$$
By considering the recurrence, we will show that $c_1=-3$ ($c_2$ can be found in a similar way).
Note that $p_2=1$, $p_3=0$ and $p_4=1/6$, moreover
$$p_n=\sum_{k=2}^{n-2} \left((k+1)p_{k+1}+{p_{k}\over k}\right)p_{n-k}{n-1 \choose k}^{-1}  \quad\mbox{for $n\geq 3$}.$$
Since $0\leq b_n\leq (n-1)!$ then $0\leq p_n\leq 1$ and
$$0\leq\sum_{k=3}^{n-5} \left((k+1)p_{k+1}+{p_{k}\over k}\right)p_{n-k}{n-1 \choose k}^{-1}\leq
\sum_{k=3}^{n-5} \left(k+2\right){n-1 \choose k}^{-1}.$$
Therefore
$$\sum_{k=3}^{n-5} \left((k+1)p_{k+1}+{p_{k}\over k}\right)p_{n-k}{n-1 \choose k}^{-1}=O\left({1\over n^3}\right)$$
because for some positive constant $C$
\begin{eqnarray*}
\sum_{k=3}^{n-5} \left(k+2\right){n-1 \choose k}^{-1}
&\leq& 5{n-1 \choose 3}^{-1}+6{n-1 \choose 4}^{-1}\\
&&+Cn^2{n-1 \choose 5}^{-1}+(n-3){n-1 \choose n-5}^{-1}
\end{eqnarray*}
Hence
\begin{eqnarray*}
{e p_n}&=&
\left(3p_{3}+{p_{2}\over 2}\right)(e p_{n-2}){n-1 \choose 2}^{-1}+O\left({1\over n^3}\right)\\
&&+\sum_{k=n-4}^{n-2} \left((k+1)(e p_{k+1})+{1\over k}(e p_{k})\right)p_{n-k}{n-1 \choose k}^{-1}\\
&=&{1\over 2}{n-1 \choose 2}^{-1}+{1\over 6}(n-3){n-1 \choose n-4}^{-1}\\
&&+\left((n-1)\left(1+{c_1\over n-1}+{c_2\over (n-1)^2}\right)+{1\over n-2}\right){n-1 \choose n-2}^{-1}
\!\!\!\!+O\left({1\over n^3}\right)\\
&=&1+{c_1\over n-1}+{c_2+3\over n^2}+O\left({1\over n^3}\right)
=1+{c_1\over n}+{c_1+c_2+3\over n^2}+O\left({1\over n^3}\right)
\end{eqnarray*}
that is
$${c_1+3\over n^2}=O\left({1\over n^3}\right)$$
which implies that $c_1=-3$.
\end{proof}

It is interesting to note that $b_n$, the number of prime expressions in $L(n)$, is related with
the number $a_n$ of {\it stabilized-interval-free} (SIF) permutations on $[n]=\{1,2,\dots,n\}$ introduced by Callan in \cite{Callan}
(seq. A075834 in \cite{Sloane}).
A permutation on $[n]$ is SIF if it does not stabilize any proper subinterval of $[n]$ . 
The SIF permutations on $[n]$ for $2\leq n\le 4$ are as follows:
$n=2\!:\,(1\,2);\ n=3\!:\,(1\,2\,3),\ (1\,3\,2)$ (the two
3-cycles); $n=4\!:\,(1\,3)(2\,4)$ and the six 4-cycles.
The power series counting the SIF permutations is

$$A(x)=\sum_{n=0}^{\infty}a_n x^n=1 + x + x^2 + 2x^3 + 7x^4 +  34x^5 + 206x^6+1476x^7+12123x^8+\cdots$$
Callan proved that
$${1\over n}[x^{n-1}](A(x))^n=(n-1)!=[x^n]F(x)$$
which means, by the Lagrange inversion formula, that 
$$A(x)={x\over F^{<-1>}(x)}=-{x\over B(x)}.$$ 
Hence
$$xA'(x)-A(x)+x=-{x\over B(x)}+{x^2B'(x)\over (B(x))^2}+{x\over B(x)}+x={x^2 B'(x)\over (B(x))^2}+x.$$
and since $xB'(x)=-(B'(x)+B(x))B(x)$ then
$$xA'(x)-A(x)+x=-{x(B'(x)+B(x))B(x)\over (B(x))^2}+x=-{xB'(x)\over B(x)}=B'(x)+B(x).$$
This differential equation yields the following recurrence formula
$$(n-1)a_n=(n+1)b_{n+1}+b_{n}\quad\mbox{for $n\geq 2$}.$$
As a consequence we confirm a numerical estimate given at the end of \cite{Callan}.

\begin{cor} \label{sif}
The asymptotic density of the SIF permutations is given by
$${a_n\over n!} = e^{-1}\left(1-{1\over n}-{5\over 2n^2}+O\left({1\over n^3}\right)\right). $$

\end{cor}

\begin{proof}

\begin{eqnarray*}
{a_n\over n!}&=&{n+1\over n-1}{p_{n+1}+{1\over n(n-1)}{p_n}}\\
&=&\left(1+{2\over n}+{2\over n^2}+O\left({1\over n^3}\right)\right){1\over e}
\left(1-{3\over n+1}-{5\over 2}{1\over (n+1)^2}+O\left({1\over n^3}\right)\right)\\
&&+{1\over n^2}{1\over e}\left(1+O\left({1\over n}\right)\right)\\
&=&{1\over e}\left(\left(1+{2\over n}+{2\over n^2}\right)
\left(1-{3\over n}\left(1-{1\over n}\right)-{5\over 2n^2}\right)+{1\over n^2}
+O\left({1\over n^3}\right)\right)\\
&=&e^{-1}\left(1-{1\over n}-{5\over 2n^2}+O\left({1\over n^3}\right)\right)
\end{eqnarray*}

\end{proof}

\end{document}